\documentclass[12pt]{amsart}
\usepackage[T1]{fontenc}
\usepackage[cp1250]{inputenc}
\usepackage[dvips]{graphicx}
\usepackage{amssymb}
\usepackage{amsmath}
\usepackage{amsthm}
\usepackage{amstext}
\usepackage{amsopn}
\usepackage{mathrsfs}
\usepackage{latexsym}
\usepackage{mathabx}
\allowdisplaybreaks
\newtheorem{Definition}{Definition}[section]
\newtheorem{Proposition}{Proposition}[section]
\newtheorem{Lemma}{Lemma}[section]
\newtheorem{Theorem}{Theorem}[section]
\newtheorem{Corollary}{Corollary}[section]
\newtheorem{Remark}{Remark}[section]
\newtheorem{Example}{Example}[section]

\numberwithin{equation}{section}
\begin{document}
%%%%%%%%%%%%%%%%%%%%%%%%%%%%%%%%%%%%%%%%%%%%%%%%%%%%%%%%%%%%
\bibliographystyle{plain}
\footnotetext{
This work is partially supported by MNiSW research grant No 1 P03A 013 30}
%%%%%%%%%%%%%%%%%%%%%%%%%%%%%%%%%%%%%%%%%%%%%%%%%%%%%%%%%%%%%%%%%
\title[Noncommutative Brownian motions associated with Kesten distributions]
{Noncommutative Brownian motions associated with Kesten distributions and related Poisson processes}
\author[R. Lenczewski]{Romuald Lenczewski}
\author[R. Sa\l{}apata]{Rafa\l{} Sa\l{}apata}
\address{Romuald Lenczewski, \newline
Instytut Matematyki i Informatyki, Politechnika Wroc\l{}awska, \newline
Wybrze\.{z}e Wyspia\'{n}skiego 27, 50-370 Wroc{\l}aw, Poland  \vspace{10pt}}
\email{Romuald.Lenczewski@pwr.wroc.pl}
\address{Rafa\l{} Sa\l{}apata, \newline
Instytut Matematyki i Informatyki, Politechnika Wroc\l{}awska, \newline
Wybrze\.{z}e Wyspia\'{n}skiego 27, 50-370 Wroc{\l}aw, Poland  \vspace{10pt}}
\email{Rafal.Salapata@pwr.wroc.pl}
%%%%%%%%%%%%%%%%%%%%%%%%%%%%%%%%%%%%%%%%%%%%%%%%%%%%%%%%%%%%%%%%
\maketitle
\begin{abstract}
We introduce and study a noncommutative two-parameter 
family of noncommutative Brownian motions in the free 
Fock space. They are associated with Kesten laws
and give a continuous interpolation between 
Brownian motions in free probability and monotone probability.
The combinatorics of our model is based
on ordered non-crossing partitions, in which to each such partition $P$ we
assign the weight $w(P)=p^{e(P)}q^{e'(P)}$,
where $e(P)$ and $e'(P)$ are, respectively, 
the numbers of {\it disorders} and {\it orders} in $P$ related
to the natural partial order on the set of blocks of 
$P$ implemented by the relation of being inner or outer. 
In particular, we obtain a simple relation between Delaney's numbers 
(related to inner blocks in non-crossing partitions) and 
generalized Euler's numbers (related to orders and disorders in 
ordered non-crossing partitions).
An important feature of our interpolation is that the mixed moments of the 
corresponding creation and annihilation processes also reproduce their monotone and 
free counterparts, which does not take place in other interpolations.
The same combinatorics is used to construct an interpolation between 
free and monotone Poisson processes. \\[5pt]

\end{abstract}
\section{Introduction}
In noncommutative probability several noncommutative Brownian motions have been 
introduced and studied. In particular, different notions of noncommutative 
independence (either abstract, or restricted to the Fock space level) 
lead to different noncommutative (i) central limit theorems, (ii) 
analogues of the classical Brownian motion, (iii) 
Poisson-type processes. In this context let us mention here 
the well known examples of the boson Brownian motion 
on the symmetric Fock space of Hudson and Parthasarathy [10],
the fermion Brownian motion on the antisymmetric Fock space of Applebaum and Hudson
[1], the free Brownian motion on the free Fock space of Speicher [22]
and the monotone Brownian motion of Muraki [18] (see also Lu [17]) on the monotone
Fock space. An important feature of the free Brownian motion is that it can be obtained as the limit in distribution of (as the dimension becomes infinite) of a sequence of Brownian motions in the finite-dimensional Hermitian matrices, as shown by Biane [4]. 

Interpolations between these examples have also been studied. 
For instance, an interpolation between the boson, 
fermion and free Brownian motions called the $q$-Brownian motion was
studied by Bo\.{z}ejko and Speicher [6]. In the bialgebra and Hopf algebra context, 
respectively, two different $q$-central limit theorems and related 
Brownian motions were studied by Schurmann [29] and the author [14,15]. 
In this paper, we introduce and study an interpolation between 
noncommutative Brownian motions in monotone probability of Muraki [19] 
and free probability of Voiculescu [22,23].

Our interpolation depends on two continuous nonnegative parameters $p,q$.
The most important reason why we find this interpolation interesting is that it is based on the 
combinatorics of ordered non-crossing partitions $\mathcal{ONC}$
and thus gives a very natural interpolation between the combinatorics 
of non-crossing partitions $\mathcal{NC}$ in free probability [21]
and the combinatorics of monotone non-crossing partitions $\mathcal{MNC}$ in monotone
probability [19]. 

By a monotone non-crossing partition of the set 
$[n]:=\{1,2, \ldots , n\}$ we understand a sequence 
$P=(P_{1},P_{2}, \ldots , P_{r})$ of blocks, in which the fact that block $P_{j}$
is inner with respect to block $P_{i}$ implies that $i<j$.
In other words, the natural partial order on the set of blocks
of $P$ implemented by the relation of being inner or outer is 
respected by the formal order defined by positions of blocks in $P$.
The term `monotone non-crossing partition' comes from our work [16], 
but it can be traced back to the partitions studied by Muraki [19]. 

The key observation is that the class of non-crossing partitions
as well as the class of monotone non-crossing partitions can be included in the 
scheme of ordered non-crossing partitions if one assigns 
the weight $p^{e(P)}$ to each $P\in  \mathcal{ONC}$, 
where $p\geq 0$ and $e(P)$ is the number of {\it disorders}, or {\it Euler inversions}.
By a disorder in $P=(P_1,P_2,\ldots , P_r)$ we understand a 
pair of blocks $\{P_{i},P_{j}\}$ such that $i<j$, $P_{i}$ is 
inner with respect to $P_{j}$ and lies immediately under $P_{j}$, i.e.
with no intermediate blocks between.
Then, for $p=0$ we obtain monotone non-crossing partitions with the same weight 
and the case $p=1$ corresponds to all ordered non-crossing partitions 
with the same weight, which reduces to $\mathcal{NC}$.

In a similar way, we can introduce the second parameter $q\geq 0$ and
assign to each $P$ the weight $q^{e'(P)}$, where $e'(P)$ is equal to the number 
of {\it orders} in $P$ (by an order we understand any 
pair $\{P_{i},P_{j}\}$ such that $i<j$ and $P_{i}$
is outer with respect to $P_{j}$). 
The second parameter is not necessary to give all Kesten laws as such
(in particular, the Wigner law) but is needed when to reproduce
all mixed moments of the free Brownian motion.
Besides, it provides a natural symmetry between orders and disorders
in our combinatorics. As a consequence, we obtain a nice relation between
{\it Delaney's numbers} $\mathcal{D}(n,k+j)$, which give 
the numbers of partitions $\pi\in \mathcal{NC}_{2n}^{2}$ 
which have $k+j$ inner blocks, and {\it generalized Euler's numbers} $\mathcal{E}(n,k,j)$, 
by which we understand the numbers of ordered partitions 
$P\in \mathcal{ONC}_{2n}^{2}$ which have $k$ disorders and $j$ orders.

Using a weight function on $L^{2}({\mathbb R}_{+})\times L^{2}({\mathbb R}_{+})$
related to $w(P)$, we define on the free Fock space $(p,q)$-creation and annihilation 
processes, $(a_{t})_{t\geq 0}$ and $(a^{*}_{t})_{t\geq 0}$, respectively, 
and the corresponding canonical position process $(\omega_{t})_{t\geq 0}$, or 
Brownian motion, where $\omega_{t} =a_{t}+a^{*}_{t}$ (parameters $p,q$ are supressed in the notation)
(in a similar way we can define the canonical momentum process $\eta_{t}=i(a_{t}-a^{*}_{t})$).
In particular, we obtain the combinatorial formula
\begin{equation}
\varphi (\omega_{t}^{2n})=
\sum_{P\in \mathcal{ONC}_{2n}^{2}}\frac{w(P)t^{n}}{b(P)!}
\end{equation}
for the even moments of the position process in the vacuum state 
$\varphi$ on ${\mathcal F}({\mathbb R}_{+})$ (the odd moments vanish), 
where the weight is given by 
\begin{equation}
w(P)=p^{e(P)}q^{e'(P)}
\end{equation}
and $b(P)$ denotes the number of blocks of $P$.
By $\mathcal{ONC}_{n}$ ($\mathcal{ONC}_{n}^{2}$) 
we denote the set of ordered non-crossing 
partitions (pair-partitions) of $[n]$.

Similarly, the moments of suitably defined processes of Poisson type
$(\gamma_{t})_{t\geq 0}$ (dependence on $p,q$ is supressed again) 
can be expressed as 
\begin{equation}
\varphi(\gamma^{n}_{t})=\sum_{P\in \mathcal{ONC}_{n}}
\frac{w(P)t^{b(P)}}{b(P)!}.
\end{equation}
It is easy to see that for $(p,q)=(0,1)$ and $(p,q)=(1,1)$ 
the above formulas give the moments of canonical position processes and
Poisson processes in monotone probability and free probability, 
respectively. The case $(p,q)=(1,0)$ corresponds to the anti-monotone processes.
Finally, the case $(p,q)=(0,0)$ corresponds to boolean processes.

Moreover, we show that our canonical position processes are associated with
Kesten distributions. For instance, for $t=1$, the moments of $\omega_{1}$ 
agree with the moments of Kesten measures with densities
\begin{equation}
f_{p,q}(x) = \frac{1}{\pi} \frac{\sqrt{2(p+q)-x^2}}{2-(2-p-q)x^2},\;\;\;x\in 
[-\sqrt{2(p+q)},\sqrt{2(p+q)}]
\end{equation}
and atoms at $x=\pm\frac{1}{\sqrt{1-(p+q)/2}}$ for $p+q<1$.
In particular, for $(p,q)=(0,1)$ (as well as for $(p,q)=(1,0)$) 
we obtain the standard arcsine law and for $(p,q)=(1,1)$ -- the standard 
Wigner law.

Let us point out that our model interpolates not only between the moments
of the canonical free and monotone position processes,
but also between the corresponding mixed moments of creation and annihilation processes.
Moreover, it reproduces independence on both Fock spaces. Therefore, at least 
on the Fock space level, it may be viewed as an interpolation between  
monotone independence and free independence. This feature is absent in the 
$t$-interpolation of Bo\.{z}ejko and Wysocza\'{n}ski [7],
which also reproduces the moments of Kesten laws [11], but
does not reproduce monotone independence. Namely, neither the mixed moments of 
monotone independent creation and annihilation processes nor the mixed moments
of position processes (with arcsine distributions) associated with disjoint
intervals agree with the corresponding moments of $t$-deformed processes
for the right value of $t$.

To put it in the general framework of interacting Fock spaces, 
let us notice that one of the main points of the $t$-interpolation and 
of the `gaussianization of probability measures' [2]
is that one can reproduce the moments of classical probability measures by means of
noncommutative Gaussian operators on the {\it one-mode interacting Fock spaces}
(see also [8,9] for a related result on symmetric measures).
Of course, one-mode interacting Fock spaces are examples of 
interacting Fock spaces [3] in which deformations of the inner product 
on the free Fock space are very simple (and are related to the Jacobi parameters 
of probability measures). Our deformations of the free Fock spaces (or, of the corresponding
creation and annihilation operators) are more complicated and, roughly speaking, 
they might be viewed as examples of {\it two-mode interacting Fock spaces}. 
Although our motivation is of combinatorial nature rather than related
to the interacting Fock space structure, it seems that this is the reason why we can 
also reproduce noncommutative independence apart from the classical 
properties like `gaussianization' in the one-mode case.

Let us mention here that a discrete interpolation between monotone probability 
and free probability, called the {\it monotone hierarchy of freeness}, was studied in [16].  
In particular, we obtained a combinatorial formula for the mixed
moments of the hierarchy of $m$-monotone Gaussians, based on the combinatorics
in which one counts blocks which are inner with respect to each block 
(as in the case of Poisson operators studied by Muraki [19]) 
rather than on the combinatorics based on blocks' depths [2], or
levels of Catalan paths [8,9], which is equivalent in the case of symmetric measures.
Moreover, this also allowed us to reproduce the mixed moments of monotone 
independent creation and annihilation operators, which are not reproduced 
by the formulas given in [2].

\section{Combinatorics of Kesten laws}

The set of all partitions of the set $[n]$ will be denoted ${\mathcal P}_n$. 
We say that $\pi\in {\mathcal P}_{n}$ is {\it non-crossing} if 
there are no pairs $\{k,k'\}\subset \pi_i$ and $\{m,m'\}\subset \pi_j$ 
with $i\neq j$ and such that $k<m<k'<m'$. 
The set of all non-crossing partitions of the set $[n]$ will be denoted $\mathcal{NC}_n$. 
By ${\mathcal P}_n^2$ we will denote the set of all {\it pair-partitions} of $[n]$ and 
$\mathcal{NC}_n^2 := \mathcal{NC}_n \cap {\mathcal P}_n^2$.

On the set of blocks of $\pi\in \mathcal{NC}_n$ we can introduce a natural partial order. Namely,
we will say that $\pi_i$ is {\it inner} with respect to $\pi_j$ for $i\neq j$ if
there exist $a,b\in\pi_j$ such that for all $c\in \pi_i$ it holds that $a < c < b$, in which case 
we shall write $\pi_j < \pi_i$. Moreover, we set $\pi_j \leqslant \pi_i $ iff $\pi_j < \pi_i$ or $\pi_j = \pi_i$.
Equivalently, we will say that $\pi_j$ is  {\it outer} w.r.t. $\pi_i$. 
We will say that blocks $\pi_j, \pi_i$ are {\it neighboring} if they are comparable 
in the above sense and there are no other `intermediate' blocks between them, 
i.e. if $\pi_j < \pi_i$ and $\pi_j\leqslant\pi_k \leqslant\pi_i$ implies that $k = i$ or $k = j$.
If $\pi_i,\pi_j$ are neighboring blocks and $\pi_j < \pi_i$, we will write $\pi_j\prec\pi_i$.
A block $\pi_i$ is called {\it outer} in $\pi$ if $\pi$ has no blocks which are outer 
with respect to $\pi_{i}$.

The pair $P=(\pi,\sigma)$, where $\pi=\{\pi_1,\pi_2,\ldots,\pi_k\}\in {\mathcal P}_n$ and 
$\sigma$ is a permutation from the symmetric group $S_k$, will be called an {\it ordered partition} 
of the set $[n]$ and will be identified with the sequence 
$P = (P_1,P_2,\ldots,P_k)$, where $P_i = \pi_{\sigma(i)}$.
In particular, we will write $P_{j}\prec P_{i}$ if $\pi_{\sigma(j)}\prec \pi_{\sigma(i)}$. 
The set of all ordered (pair, non-crossing, non-crossing pair)
partitions of $[n]$ will be denoted ${\mathcal OP}_n$ (respectively, 
${\mathcal OP}_n^2$, ${\mathcal{ONC}}_n$, $\mathcal{ONC}_n^2$).

Let us observe that in each ordered partition $P=(\pi,\sigma)$, the permutation 
$\sigma$ defines a linear order on the set of blocks of $\pi$.
Comparing  this order with the partial order given by the relation of being inner (outer)
for $P\in\mathcal{ONC}_n$, we can introduce ,,disorders'' betwen blocks.

\begin{Definition}\label{d18}
{\rm
If $P=(P_1,P_2,\ldots,P_k)\in\mathcal{ONC}_n$, we will say that the pair $\{P_i,P_j\}$ forms a
{\it disorder} or {\it Euler inversion} if $i < j$ and $P_j \prec P_i$. 
If $i<j$ and $P_{i}\prec P_{j}$, we will say that the pair $\{P_{i},P_{j}\}$ 
forms an {\it order}.
The numbers of all disorders and orders in $P$ will be denoted $e(P)$ and $e'(P)$,
respectively.
In a similar way we define disorders and orders in the permutation $\sigma\in S_n$ 
associated with each index $i\in\{1,2,\ldots,n-1\}$ for which  
$\sigma(i) < \sigma(i+1)$ and $\sigma(i)>\sigma(i+1)$, respectively.
The numbers of all disorders and orders in $\sigma$ will be denoted $e(\sigma)$ and $e'(P)$,
respectively.}
\end{Definition}
\begin{Remark}
{\rm Well-known Euler numbers, denoted $\left< \genfrac{}{}{0pt}{}{n}{k} \right>$,
give the numbers of permutations of the set $[n]$ which have $k$ disorders.}
\end{Remark}
\begin{figure}
\unitlength=0.8mm
\special{em:linewidth 0.4pt}
\linethickness{0.4pt}
\begin{picture}(111.00,45.00)(10.00,20.00)

\put(0.00,41.00){\line(0,1){14.00}}
\put(0.00,55.00){\line(1,0){42.00}}
\put(42.00,41.00){\line(0,1){14.00}}
\put(20.00,56.00){\footnotesize{\it{4}}}

\put(6.00,41.00){\line(0,1){4.00}}
\put(6.00,45.00){\line(1,0){6.00}}
\put(12.00,41.00){\line(0,1){4.00}}
\put(8.50,46.00){\footnotesize{\it{1}}}

\put(18.00,41.00){\line(0,1){9.00}}
\put(18.00,50.00){\line(1,0){18.00}}
\put(36.00,41.00){\line(0,1){9.00}}
\put(26.00,51.00){\footnotesize{\it{2}}}

\put(24.00,41.00){\line(0,1){4.00}}
\put(24.00,45.00){\line(1,0){6.00}}
\put(30.00,41.00){\line(0,1){4.00}}
\put(26.00,46.00){\footnotesize{\it{3}}}

\put(0.00,41.00){\circle*{1.00}}
\put(6.00,41.00){\circle*{1.00}}
\put(12.00,41.00){\circle*{1.00}}
\put(18.00,41.00){\circle*{1.00}}
\put(24.00,41.00){\circle*{1.00}}
\put(30.00,41.00){\circle*{1.00}}
\put(36.00,41.00){\circle*{1.00}}
\put(42.00,41.00){\circle*{1.00}}

\put(-0.90,37.00) {\footnotesize 1}
\put(05.10,37.00) {\footnotesize 2}
\put(11.10,37.00) {\footnotesize 3}
\put(17.10,37.00) {\footnotesize 4}
\put(23.10,37.00) {\footnotesize 5}
\put(29.10,37.00) {\footnotesize 6}
\put(35.10,37.00) {\footnotesize 7}
\put(41.10,37.00) {\footnotesize 8}

\put(96.00,41.00){\line(0,1){4.00}}
\put(90.00,55.00){\line(1,0){42.00}}
\put(132.00,41.00){\line(0,1){14.00}}
\put(110.00,56.00){\footnotesize{\it{1}}}

\put(90.00,41.00){\line(0,1){14.00}}
\put(108.00,41.00){\line(0,1){9.00}}
\put(95.00,46.00){\footnotesize{\it{3}}}

\put(102.00,41.00){\line(0,1){9.00}}
\put(102.00,50.00){\line(1,0){24.00}}
\put(126.00,41.00){\line(0,1){9.00}}
\put(113.00,51.00){\footnotesize{\it{2}}}

\put(114.00,41.00){\line(0,1){4.00}}
\put(114.00,45.00){\line(1,0){6.00}}
\put(120.00,41.00){\line(0,1){4.00}}
\put(116.00,46.00){\footnotesize{\it{4}}}

\put(90.00,41.00){\circle*{1.00}}
\put(96.00,41.00){\circle*{1.00}}
\put(102.00,41.00){\circle*{1.00}}
\put(108.00,41.00){\circle*{1.00}}
\put(114.00,41.00){\circle*{1.00}}
\put(120.00,41.00){\circle*{1.00}}
\put(126.00,41.00){\circle*{1.00}}
\put(132.00,41.00){\circle*{1.00}}

\put(89.10,37.00) {\footnotesize 1}
\put(95.10,37.00) {\footnotesize 2}
\put(101.10,37.00) {\footnotesize 3}
\put(107.10,37.00) {\footnotesize 4}
\put(113.10,37.00) {\footnotesize 5}
\put(119.10,37.00) {\footnotesize 6}
\put(125.10,37.00) {\footnotesize 7}
\put(131.10,37.00) {\footnotesize 8}
\put(78,28){\small{$Q=\big(\{1,8\},\{3,4,7\},\{2\},\{5,6\}\big)$}}
\put(-12,28){\small{$P=\big(\{2,3\},\{4,7\},\{5,6\},\{1,8\}\big)$}}
\end{picture}
\caption{Examples of ordered non-crossing partitions}
\end{figure}
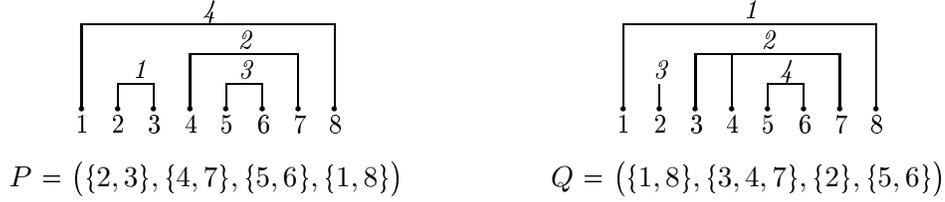

\begin{Example}\label{p9}
{\rm
Consider ordered non-crossing partitions 
$P\in\mathcal{ONC}_8^2$ and $Q\in\mathcal{ONC}_8$ given in Fig.1.
Parititon $P$ has 3 pairs of neighboring blocks: 
$P_4 \prec P_2$, $P_4 \prec P_1$ and $P_2 \prec P_3$. 
Only the first two give disorders, thus $e(P) = 2$ and $e'(P)=1$.
Pairs of neighboring blocks in $Q$ are the following: $Q_1 \prec Q_2$, $Q_1 \prec Q_3$ 
and $Q_2 \prec Q_4$. None of them gives a disorder, thus $e(Q) = 0$ and $e'(Q)=3$.
Non-crossing ordered paritions which do not have disorders are called
{\it monotone} [16].
}
\end{Example}

By ${\mathcal {ONCC}}_n\ ({\mathcal {ONCC}}_n^2)$ we denote the set of those ordered 
non-crossing partitions (pair-partitions) of $[n]$, in which
the numbers $1$ and $n$ belong to the same block. Such partitions will be called 
{\it covered}. Let us introduce numbers
$$
r_n = \frac{1}{n!} \sum_{P\in \mathcal{ONC}_{2n}^2} w(P)\ ,\ \ \ \ \ \ 
s_n = \frac{1}{n!} \sum_{P\in \mathcal{ONCC}_{2n}^2} w(P) 
$$
for $n\geq 1$ and set $r_0 =1$, $s_0 = 0$, where $w(P)$ is given by (1.2).

\begin{Proposition}
The following relation between sequences $(r_n)_{n=1}^{\infty}$ and
$(s_n)_{n=1}^{\infty}$ holds:
$$
r_n =\sum_{m=1}^n \sum_{k_1 + k_2 + \ldots + k_m = n}  s_{k_1} s_{k_2}\,\ldots\, s_{k_m}\ ,\ \ \ \ \ \ \ \ n\geqslant1.
$$
where the second summation runs over positive indices $k_1,\ldots , k_m$.
\end{Proposition}

\noindent
{\bf Proof.} Let $Q^{(i)} = (Q^{(i)}_1,Q^{(i)}_2,\ldots,Q^{(i)}_{k_i})$, where $i=1,2,\ldots,m$,
be arbitrary pair-partitions from the sets $\mathcal{ONCC}_{2k_i}^2$, respectively, such that
$k_1 + \ldots  + k_m = n$. From the blocks of all these partitions we construct 
ordered pair partitions of $[n]$ with the shape given by Fig.2. We order all blocks $Q^{(i)}_j$ 
of the subpartition $Q^{(i)}$ in such a way that the order 
between blocks from the same partition $Q^{(i)}$ is preserved. 
There are $\frac{n!}{k_1!\ldots k_m!}$ such orderings and each of them
defines exactly one $P\in \mathcal{ONC}_{2n}^2$. 
Moreover, each partition from $\mathcal{ONC}_{2n}^2$ can be obtained in this fashion
by an appropriate choice of covered partitions $Q^{(i)}$. 
From the above reasoning we obtain
\begin{equation}\label{r30}
\left| \mathcal{ONC}_{2n}^2 \right| = \sum_{m=1}^n \sum_{k_1 + k_2+\ldots + k_m = n} 
\frac{n!}{k_1!\ldots k_m!} 
\left| \mathcal{ONCC}_{2k_1}^2 \right| 
\ldots 
\left| \mathcal{ONCC}_{2k_m}^2 \right|.
\end{equation}

\begin{figure}
\unitlength=1mm
\special{em:linewidth 0.4pt}
\linethickness{0.4pt}
\begin{picture}(120.00,40.00)(-30.00,15.00)\label{rys2}

\put(-19.70,46.00){$\overbrace{\hspace{282pt}}^{P}$}

\put(-18.00,31.00){\line(0,1){12.00}}
\put(-18.00,43.00){\line(1,0){25.00}}
\put(7.00,31.00){\line(0,1){12.00}}
\put(-18.00,31.00){\circle*{1.00}}
\put(7.00,31.00){\circle*{1.00}}
\put(-7.00,31.00){$\cdots$}
\put(-19.70,27.00){$\underbrace{\hspace{81pt}}_{Q^{(1)}}$}

\put(13.00,31.00){\line(0,1){12.00}}
\put(13.00,43.00){\line(1,0){25.00}}
\put(38.00,31.00){\line(0,1){12.00}}
\put(13.00,31.00){\circle*{1.00}}
\put(38.00,31.00){\circle*{1.00}}
\put(24.00,31.00){$\cdots$}
\put(11.30,27.00){$\underbrace{\hspace{81pt}}_{Q^{(2)}}$}

\put(44.00,31.00){$\cdots$}

\put(53.00,31.00){\line(0,1){12.00}}
\put(53.00,43.00){\line(1,0){25.00}}
\put(78.00,31.00){\line(0,1){12.00}}
\put(53.00,31.00){\circle*{1.00}}
\put(78.00,31.00){\circle*{1.00}}
\put(64.00,31.00){$\cdots$}
\put(51.30,27.00){$\underbrace{\hspace{81pt}}_{Q^{(m)}}$}

\end{picture}
\caption{$P\in \mathcal{ONC}$ constructed from $Q^{(1)}, \ldots , Q^{(m)}
\in \mathcal{ONCC}$.}
\end{figure}
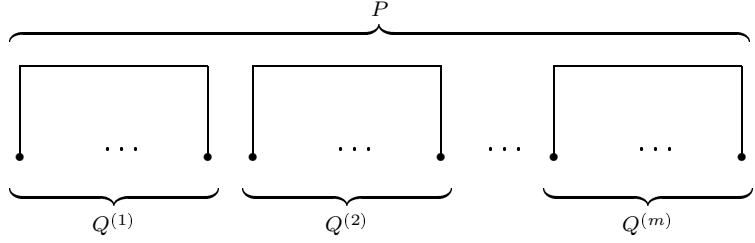
Clearly, between blocks which belong to different 
partitions $Q^{(i)}$ and $Q^{(j)}$, $i\neq j$, there are no orders or disorders.
Therefore,
\begin{eqnarray*}
e(P) &=& e(Q^{(1)}) + e(Q^{(2)}) + \ldots + e(Q^{(r)})\\
e'(P) &=& e'(Q^{(1)}) + e'(Q^{(2)}) + \ldots + e'(Q^{(r)})
\end{eqnarray*}
which gives multiplicativity of the weights
$$
w(P)=w(Q^{(1)})w(Q^{(2)})\ldots w(Q^{(r)})
$$
and that, together with (\ref{r30}), gives the assertion.$\hfill{\blacksquare}$

\begin{Corollary}\label{w3}
Let 
$R(z) = \sum_{n=0}^{\infty} r_n \,z^n$ and 
$S(z) = \sum_{n=0}^{\infty} s_n \,z^n\,$ be formal power series.
Then
$$
R(z) = \frac{1}{1-S(z)}\,.
$$
\end{Corollary}

\noindent
{\bf Proof.} Using Proposition 2.1 and simple algebraic computations, we get
\begin{eqnarray*}
R(z)-1 &=& \sum_{n=1}^{\infty} r_n z^n = 
\sum_{n=1}^{\infty} \sum_{m=1}^n 
\sum_{k_1 + k_2 + \ldots + k_m = n}  
s_{k_1} s_{k_2}\,\ldots\, s_{k_m} z^n \\
&=& 
\sum_{m=1}^{\infty} \sum_{n=m}^{\infty} 
\sum_{k_1 + k_2 + \ldots + k_m= n}  
s_{k_1} s_{k_2}\,\ldots\, s_{k_m} z^n \\
&=& 
\sum_{m=1}^{\infty}(S(z))^{m}= \frac{S(z)}{1-S(z)}
\end{eqnarray*}
from which our assertion follows. \hfill $\blacksquare$\\

\begin{figure}
\unitlength=1mm
\special{em:linewidth 0.4pt}
\linethickness{0.4pt}
\begin{picture}(120.00,40.00)(-30.00,15.00)\label{rys1}

%\put(-19.70,46.00){$\overbrace{\hspace{282pt}}^{\pi}$}

\put(-24.00,31.00){\circle*{1.00}}
\put(-24.00,31.00){\line(0,1){16.00}}
\put(-24.00,47.00){\line(1,0){108.00}}
\put(84.00,31.00){\circle*{1.00}}
\put(84.00,31.00){\line(0,1){16.00}}
\put(29.00,49.00){${\footnotesize P_{n+1}}$}

\put(-18.00,31.00){\line(0,1){12.00}}
\put(-18.00,43.00){\line(1,0){25.00}}
\put(7.00,31.00){\line(0,1){12.00}}
\put(-18.00,31.00){\circle*{1.00}}
\put(7.00,31.00){\circle*{1.00}}
\put(-7.00,31.00){$\cdots$}
\put(-19.70,27.00){$\underbrace{\hspace{81pt}}_{Q_1}$}

\put(13.00,31.00){\line(0,1){12.00}}
\put(13.00,43.00){\line(1,0){25.00}}
\put(38.00,31.00){\line(0,1){12.00}}
\put(13.00,31.00){\circle*{1.00}}
\put(38.00,31.00){\circle*{1.00}}
\put(24.00,31.00){$\cdots$}
\put(11.30,27.00){$\underbrace{\hspace{81pt}}_{Q_2}$}

\put(44.00,31.00){$\cdots$}

\put(53.00,31.00){\line(0,1){12.00}}
\put(53.00,43.00){\line(1,0){25.00}}
\put(78.00,31.00){\line(0,1){12.00}}
\put(53.00,31.00){\circle*{1.00}}
\put(78.00,31.00){\circle*{1.00}}
\put(64.00,31.00){$\cdots$}
\put(51.30,27.00){$\underbrace{\hspace{81pt}}_{Q_r}$}

\end{picture}
\caption{Partition covered by the last block $P_{n+1}$.}
\end{figure}
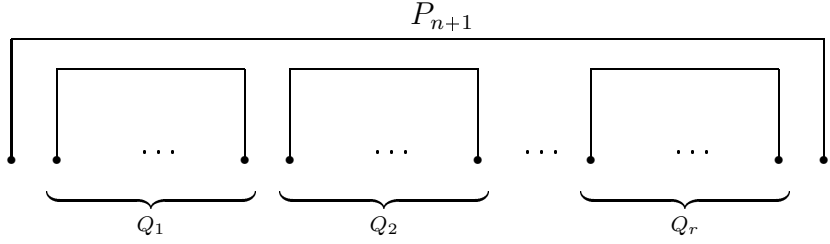

In order to find $R(z)$, we introduce
another sequence of numbers, denoted $(a_n)_{n=0}^{\infty}$ and defined by the combinatorial formula
\begin{equation}\label{r31}
a_n = \sum_{m=1}^{n} \sum_{k_1 + k_2 + \ldots + k_m = n} p^m  
s_{k_1} s_{k_2}\,\ldots\, s_{k_m}
\end{equation}
for $n\geq 1$ and we set $a_0=1$. 
Let us observe that $a_n$ is the sum of contributions to $r_{n+1}$ of these
pair partitions $P=(P_1,\ldots ,P_{n+1})$ $\in \mathcal{ONCC}_{2n+2}^2$
which are covered by the block of highest color, namely $P_{n+1}$. 
Therefore,
\begin{equation}\label{r32}
a_n =
\ \frac{1}{n!} 
\sum_{\stackrel{P\in \mathcal{ONCC}_{2n+2}^2}{\scriptscriptstyle P_{n+1} = \{1,2n+2\}}} w(P)
\end{equation}
for $n\geq 0$.
In fact, each neighboring block of $P_{n+1}$ corresponds to
a certain covered partition $Q_i\in \mathcal{ONCC}_{2k_i}^2$ for 
$i=1,\ldots ,m$, as Fig.3 shows. 
Of course, $k_1+\ldots+k_m = n$. 
Besides, let us observe that each neighboring block of $P_{n+1}$ forms a disorder
with block $P_{n+1}$.
This justifies equivalence of (\ref{r31}) and (\ref{r32}).
\begin{Corollary}\label{w4}
Let $A(z) = \sum_{n=0}^{\infty} a_n \, z^n$ be a formal power series. Then
$$
A(z) = \frac{1}{1-pS(z)} \,.
$$
\end{Corollary}

\noindent
{\bf Proof.} 
Using algebraic calculations and equation (\ref{r31}), we obtain
\begin{eqnarray*}
A(z)-1 &=& \sum_{n=1}^{\infty} a_n z^n = 
\sum_{n=1}^{\infty} \sum_{m=1}^n 
\sum_{k_1 + k_2 + \ldots + k_m = n}
p^m s_{k_1} s_{k_2}\,\ldots\, s_{k_m} z^n \\
&=& 
\sum_{m=1}^{\infty} p^m 
\sum_{n=m}^{\infty} 
\sum_{k_1 + k_2 + \ldots + k_m = n} s_{k_1} s_{k_2}\,\ldots\, s_{k_m} z^n \\
&=& 
\sum_{m=1}^{\infty} \left( pS(z) \right)^m 
= \frac{pS(z)}{1-pS(z)}
\end{eqnarray*}
from which we get the assertion.\hfill $\blacksquare$\\

To find a relation between $R(z)$ and $A(z)$ we shall need some additional notations.
Let $\mathcal{ONCC}_{n}(r)$ ($\mathcal{ONCC}^2_{n}(r)$) denote the sets of
ordered partitions (pair partitions) of $[n]$ with $r$ outer blocks, and, for  $r=1,\ldots,n$,
introduce sequences 
\begin{equation}\label{r1}
s_n^{(r)} = \frac{1}{n!} \sum_{Q\in {\mathcal{ONCC}}_{2n}^2(r)} w(P)
\end{equation}
for $n\geq 1$, and set $s_0^{(r)}=0$. Of course, $s_n^{(1)}=s_n$ and $s^{(r)}_r = 1$.

\begin{Proposition}\label{f1}
Let $S^{(r)}(z) = \sum_{n=r}^{\infty} s_n^{(r)} z^n$, where $r\geq 1$. Then 
$$
S^{(r)}(z) = \big( S(z) \big)^r
$$
\end{Proposition}

\noindent
{\bf Proof.} Notice that we have the following relation
between sequences $(s_n)_{n\geq 1}$ and $(s_n^{(r)})_{n\geq 1}$:
$$
s_n^{(r)} = \sum_{k_1 + k_2 + \ldots + k_r = n} s_{k_1} \ldots s_{k_r}\,,\ \ \ \ \ \ \ \ n\geqslant r\,,
$$
where $k_{1},k_{2}, \ldots, k_{r}$ are assumed to be positive integers.
This gives
$$
S^{(r)}(z) = 
\sum_{n=r}^{\infty} \sum_{k_1 + k_2 + \ldots + k_r = n} 
s_{k_1} \ldots s_{k_r} z^n
= \Big( \sum_{n=1}^{\infty} s_n z^n \Big)^r = \big( S(z) \big)^r\,.
$$
which proves our assertion.
$\hfill{\blacksquare}$

\begin{Lemma}\label{l1}
For $n\geqslant r+1$ it holds that
$$
s_n^{(r)} = \frac{r}{n} \sum_{k=1}^{n-r+1} a_{k-1} s_{n-k}^{(r-1)} \ +\  \frac{q}{n} \sum_{k=1}^{n-r} (2n-2k-r) a_{k-1} s_{n-k}^{(r)}\,.
$$
\end{Lemma}

\noindent
{\bf Proof.} 
Let us split (\ref{r1}) into two sums: the first one running over
those partitions from ${\mathcal {ONCC}}_{2n}^2(r)$ in which the block of highest color
is outer, and the second one -- over the remaining partitions. Then
\begin{eqnarray*}
s_n^{(r)} &=& \frac{1}{n!} \sum_{\stackrel{P\in{\mathcal{ONCC}}_{2n}^2(r)}{\scriptscriptstyle P_n - outer}} w(P) \ +\ \frac{1}{n!} \sum_{\stackrel{P\in{\mathcal{ONCC}}_{2n}^2(r)}{\scriptscriptstyle P_n - inner}} 
w(P) \\
&=& \frac{1}{n!} \sum_{k=1}^{n-r+1} \sum_{\stackrel{T\in{\mathcal{ONCC}}_{2k}^2}{\scriptscriptstyle T_k = \{1,2k\}}} \sum_{\scriptscriptstyle Q\in{\mathcal{ONCC}}_{2n-2k}^2(r-1)} r \left( \genfrac{}{}{0pt}{}{n-1}{k-1} \right) w(T)w(Q) \\
&+& \frac{1}{n!} \sum_{k=1}^{n-r} \sum_{\stackrel{T\in{\mathcal{ONCC}}_{2k}^2}{\scriptscriptstyle T_k = \{1,2k\}}} \sum_{\scriptscriptstyle Q\in{\mathcal{ONCC}}_{2n-2k}^2(r)} (2n-2k-r) \left( \genfrac{}{}{0pt}{}{n-1}{k-1} \right) qw(T)w(Q) \\
&=& \frac{r}{n} \sum_{k=1}^{n-r+1} a_{k-1} s_{n-k}^{(r-1)} \ +\  \frac{q}{n} \sum_{k=1}^{n-r} (2n-2k-r) a_{k-1} s_{n-k}^{(r)}\,.
\end{eqnarray*}
which completes the proof.$\hfill{\blacksquare}$

\vspace{12pt}
From the above lemma we get a relation between the $S^{(r)}(z)$'s and $A(z)$
in the form of a differential equation. 

\begin{Corollary}\label{w5}
The functions $S^{(r)}(z)$ and $A(z)$ satisfy the differential
recurrence 
$$
\big( S^{(r)}(z)\big)' = r S^{(r-1)}(z) A(z) + 2qz \big( S^{(r)}(z)\big)' A(z) - qrA(z)S^{(r)}(z)\,,
$$
with initial conditions $S^{(r)}(0) = 0\,, A(0) = 1$.
\end{Corollary}

\noindent
{\bf Proof.} In view of Lemma \ref{l1}, we have
\begin{eqnarray*}
\big( S^{(r)}(z)\big)' &=& \sum_{n=r}^{\infty} n s^{(r)}_{n} z^{n-1} \\
&=& r z^{r-1} + r \sum_{n=r+1}^{\infty} \sum_{k=1}^{n-r+1} a_{k-1} s_{n-k}^{(r-1)} z^{n-1} \\
& & + q\sum_{n=r+1}^{\infty} \sum_{k=1}^{n-r} (2n-2k-r) a_{k-1} s^{(r)}_{n-k} z^{n-1} \\
&=& r \sum_{n=r}^{\infty} \sum_{k=1}^{n-r+1} a_{k-1} s_{n-k}^{(r-1)} z^{n-1} + 2q \sum_{n=r+1}^{\infty} \sum_{k=1}^{n-r} (n-k) a_{k-1} s^{(r)}_{n-k} z^{n-1} \\
& & - q r \sum_{n=r+1}^{\infty} \sum_{k=1}^{n-r} a_{k-1} s^{(r)}_{n-k} z^{n-1} \\
&=& r \sum_{n=r-1}^{\infty} \sum_{k=0}^{n-(r-1)} a_{k} s_{n-k}^{(r-1)} z^{n} + 2qz \sum_{n=r}^{\infty} \sum_{k=0}^{n-r} a_{k} (n-k) s^{(r)}_{n-k} z^{n-1} \\
& & - q r \sum_{n=r}^{\infty} \sum_{k=0}^{n-r} a_{k} s^{(r)}_{n-k} z^{n}\\
&=& r S^{(r-1)}(z) A(z) + 2qz \big( S^{(r)}(z)\big)' A(z) - qrA(z)S^{(r)}(z)
\end{eqnarray*}
Of course, $S^{(r)}(0)=0$ and $A(0)=a_0=1$, which completes the proof. $\hfill{\blacksquare}$\\
\indent{\par}
Using the relations bewteen functions $R(z), A(z)$ and $S^{(r)}(z)$, $r\geq 1$, 
we can now derive the explicit form of $R(z)$, which will turn out to be related 
to the moment generating function of Kesten laws.

\begin{Theorem}\label{t15}
For each $p\geq 0, p\geq 0$, $p+q>0$, the sequence $(m_{n})_{n\geq 0}$,
where 
$$
m_{n}=\left\{
\begin{array}{ll}
r_{k} & if\;\;n=2k\\
0     & if\;\;n\;\;is\;\; odd
\end{array}
\right.
$$
is the sequence of moments of the Kesten measure (1.4). 
The probability measure determined by the moments is unique.
\end{Theorem}

\noindent
{\bf Proof.} Using Corollary \ref{w5} and Proposition 2.2, 
we obtain a differential equation for $S(z)$, namely
$$
rS'(z)=rA(z)+2qrzS'(z)A(z)-qrS(z)A(z)
$$
which, in view of Corollary \ref{w3}, leads to the 
differential equation for $R(z)$ of the form
$$
R'(z) = \frac{R^2(z)\big( (1-q)R(z) +q \big)}{R(z) \big( 1-p-2qz \big)+p}
$$
with the initial condition $R(0)=1$.
Let us observe now that the function $R(z)$ must be symmetric with respect to $p$ and $q$
(this easily follows from the definition of $R(z)$ if we reverse the order in all
ordered non-crossing pair partitions over which the summation is taken).
Thus, if we denote $R(z)=R_{p,q}(z)$, then $R_{p,q}(z)=R_{q,p}(z)$.
Therefore, we obtain another differential equation for $R(z)$ with $p$ and $q$ interchanged.
This allows us then to reduce the above
equation to the (algebraic) quadratic equation, namely
$$
AR^{2}(z)+BR(z)+C=0
$$
where 
\begin{eqnarray*}
A&=&(1-q)^2-(1-p)^2-2pz+2qz\,,\\
B&=&2q(1-q)-2p(1-p)\,,\\
C&=&q^2-p^2\,,
\end{eqnarray*}
which has two solutions, 
$$
R(z) = \frac{p+q-1 \pm \sqrt{1-2(p+q)z}}{p+q-2+2z}\,
$$
but only the one corresponding to the minus sign satisfies the initial condition 
$R(0)=1$.
If we define the corresponding moment generating function by taking
$m_{n}=r_{k}$ for $n=2k$ with odd moments equal to zero, we obtain 
the function $M(z)=\sum_{n=0}^{\infty}m_{n}z^{n}=R(z^2)$ and the corresponding
Cauchy transform 
$$
G(z) = \frac{1}{z} M\left(\frac{1}{z}\right) = \frac{(p+q-1)z-\sqrt{z^2-2(p+q)}}{2-(2-p-q)z^2}\,.
$$ 
turns out to be the Cauchy transform of the (uniquely determined) 
Kesten distribution $\mu_{p,q}$ (with the absolutely continuous part given by (1.4)). 
$\hfill{\blacksquare}$

\begin{Remark}\label{f8}
{\rm The continued fraction representation of $G(z)$ takes the form
$$
G(z) = \cfrac{1}{z-\cfrac{1}{z- \cfrac{t}{z - \cfrac{t}{z - \ddots}}}}\ .
$$
where $t=1/2(p+q)$. Such Cauchy transforms were obtained 
in the context of the so-called $t$-transformation of measures [7]. However,
we will show later that our combinatorics is different and, when carried over
to the Fock space level, gives a different Brownian motion (although it also
has Kesten distributions).}
\end{Remark}

\section{Noncommutative Brownian motions}\label{s6}

We will now construct new types of noncommutative 
Brownian motions on the free Fock space which have Kesten distributions and are parametrized 
by two nonnegative real numbers $p,q$. They can be viewed as an interpolation between
the free Brownian motion obtained for $(p,q)=(1,1)$
and the monotone Brownian motion corresponding to 
$(p,q)=(0,1)$ (the anti-monotone and boolean Brownian motions are also obtained,
for $(p,q)=(1,0)$ and $(p,q)=(0,0)$, respectively).
Moreover, the mixed moments of the associated creation and 
annihilation operators in the vacuum state agree with their counterparts
in free probability and monotone probability, a feature absent in other interpolations.

By the free Fock space over a Hilbert space ${\mathcal H}=L^{2}({\mathbb R}_{+})$
we understand the direct sum
\begin{equation}\label{r21}
{\mathcal F}({\mathcal H})=
{\mathbb C}\Omega \oplus  \bigoplus_{n=1}^{\infty}{\mathcal H}^{\otimes n}  \cong {\mathbb C}\Omega \oplus \bigoplus_{n=1}^{\infty}L^{2}({\mathbb R}_{+}^{n})
\end{equation}
with the canonical inner product.

Let $w$ be the weight function on ${\mathbb R}_{+}\times {\mathbb R}_{+}$ given by 
$$
w(s,t) = \left\{
\begin{array}{ll}p & {\rm if}\;\; 0<s<t\\
q& {\rm if} \;\; 0<t< s\\
1 & {\rm otherwise}
\end{array}
\right. ,
$$
where $p>0$, $q\geq 0$.
Now, introduce special vectors in ${\mathcal F}({\mathcal H})$
denoted by the `tensor-like' symbol $f_1\oasterisk f_2\oasterisk \ldots \oasterisk f_n$, where
$$
(f_1\oasterisk f_2 \oasterisk \ldots \oasterisk f_n)
(t_1,t_2, \ldots , t_n):=
f_1(t_1)\sqrt{w(t_1,t_2)}f_2(t_2)\ldots \sqrt{w(t_{n-1},t_n)}f_n(t_n).
$$
These vectors remind simple tensors, but they have the 
`nearest neighbor coupling'.
Note that the set of such vectors is dense in ${\mathcal F}({\mathcal H})$
if $p>0$ and $q>0$. If $(p,q)=(0,1)$, 
it is dense in the monotone Fock space
$$
{\mathcal M}({\mathcal H}) = 
{\mathbb C}\Omega \oplus \bigoplus_{n=1}^{\infty} L^{2} (\Delta^{(n)})
$$
where $\Delta^{(n)} = \{ (t_1,\ldots,t_n)\in\mathbb{R}_{+}^n; t_{1} \leqslant \ldots \leqslant t_{n} \}$.
In turn, if $(p,q)=(1,0)$, it is dense in the anti-monotone Fock space 
(similar to the monotone Fock space, but with the reversed order of coordinates).
If $(p,q)=(0,0)$, we obtain in turn 
${\mathbb C}\Omega \oplus L^{2}({\mathbb R}_{+})$.

Let us define suitable creation, gauge and annihilation operators 
which are $(p,q)$-deformations of their free counterparts.

\begin{Definition}
{\rm Define the $(p,q)$-{\it creation operator}
$a(f): {\mathcal F}({\mathcal H}) \rightarrow {\mathcal F}({\mathcal H})$
associated with $f\in{\mathcal H}$ as the bounded linear extension of
\begin{eqnarray}
a(f)\,\Omega &=& f\\
a(f)\, \left(f_1\oasterisk f_2 \oasterisk \ldots  \oasterisk f_n \right)
&=&  
f\oasterisk f_1\oasterisk \ldots \oasterisk f_n
\end{eqnarray}
where  $f,f_1, \ldots , f_n\in {\mathcal H}$.}
\end{Definition}

\begin{Definition}
{\rm For any $f\in L^{\infty}({\mathbb R}_{+})$ with $|f(0)|<\infty$, by  
the $(p,q)$-{\it gauge operator} on ${\mathcal F}({\mathcal H})$ 
associated with $f$ we understand the bounded linear extension of 
\begin{eqnarray*}
M(f)\,\Omega&=&f(0)\Omega\\
M(f)\,\left(f_1 \oasterisk f_2\oasterisk \ldots \oasterisk f_n\right)&=
&(ff_1)\oasterisk f_2\oasterisk \ldots \oasterisk f_n
\end{eqnarray*}
for any $f_1,f_2, \ldots , f_n \in {\mathcal H}$.}
\end{Definition} 

In particular, by $M(f,g)$ we will denote the gauge operator
associated with the function defined by the weighted inner product on 
${\mathcal H}$, namely
$$
\langle \!\langle f,g\rangle \!\rangle (t):=
\int_{{\mathbb R}_{+}}f(s)\overline{g(s)}w (s,t)ds \,.
$$
Note that this gauge operator multiplies the vacuum vector by the value of the 
inner product of $f$ and $g$, i.e. 
\begin{equation}
M(f,g)\,\Omega=\langle f, g \rangle \,\Omega
\end{equation}
thanks to our assumption that $w(s,0)=1$ for any $s\geq 0$. 
The $(p,q)$-gauge operator allows us to write a simple formula for the
action of the annihilation operators.

\begin{Proposition}\label{f10}
The action of the $(p,q)$-{\it annihilation operator} 
$a^*(f)$ associated with $f\in {\mathcal H}$, 
adjoint with respect to the creation operator
$a(f)$, is given by the bounded linear extension of
\begin{eqnarray}
a^*(f)\, \Omega 
&=& 
0 \nonumber \\
a^*(f) \,\left(f_1\oasterisk f_2\oasterisk \ldots \oasterisk f_n \right)
&=& 
M(f_1,f)\,\left(f_2\oasterisk f_3\oasterisk \ldots \oasterisk f_{n}\,,\label{r34}\right)
\end{eqnarray}
where $f_1,f_2, \ldots , f_n\in {\mathcal H}$.
Thus, in particular, $a^*(f) f_1 = \langle f_1, f \rangle$.\\ 
\end{Proposition}

\begin{Remark}
{\rm Equivalently, we can use a deformed inner product on the free Fock space
and define the creation operators to be the free creation operators, whereas 
the annihilation operators to be their adjoints with respect to the $(p,q)$-deformed 
inner product given by 
$$
\langle F , G \rangle 
= 
\delta_{n,m} \sum_{\sigma\in S_n} 
w(\sigma^{-1}) 
\int_{\Delta_{\sigma}} F(t_1,\ldots,t_n)\overline{G(t_1,\ldots,t_n)} dt_1 \ldots dt_n\,,
$$
for any $p,q>0$, where $F\in L^{2}({\mathbb R}_{+}^n)$, $G\in L^{2}({\mathbb R}_{+}^{m})$,
and $\langle \Omega, \Omega \rangle=1$, $\langle \Omega, F, \rangle =0$.
This definition can be extended to $p,q\geq 0$ except that one has
to divide the above vector space 
by the corresponding kernel of the sesquilinear form.}
\end{Remark}

The gauge operators (which commute among themselves) allow us to write relations between creation, 
annihilation and gauge operators in a simple form as the proposition 
given below shows (we omit the elementary proof).

\begin{Proposition}
The following relations hold:
\begin{equation}
a^{*}(g)a(f)=M(f,g), \;\;\;
M(h)a(f)=a(hf),\;\;\;
a^{*}(f)M(\overline{h})=a^{*}(fh)
\end{equation}
where $f,g\in {\mathcal H}$ and $h\in L^{\infty}({\mathbb R}_{+})$ 
with $|h(0)|<\infty$.
\end{Proposition}

By the $(p,q)$-{\it Gaussian operator} associated with the function $f\in {\mathcal H}$ 
we will understand the self-adjoint position operator given by the sum $\,\omega(f) = a^*(f) + a(f)$. 
In turn, by the associated {\it Brownian motion} we will understand 
the process $(\omega_{t})_{t\geq 0}$, where $\omega_{t}=\omega(\chi_{[0,t)})$.
Below we will find a formula for the mixed moments
$\varphi(\omega(f_1) \omega(f_2) \ldots \omega(f_n))$,
where $\varphi$ is the vacuum state on ${\mathcal F}({\mathcal H})$ 
and $ f_1,f_2,\ldots,f_n \in \Theta$, where
\begin{equation}
\Theta:=\{ \chi_{[s,t)}; 0 \leqslant s < t < \infty \}
\end{equation}
is the set of characteristic functions of intervals.

We shall assume that supports of these functions are pairwise disjoint and are ordered by the 
partial order $I_1<I_2$ whenever $t_1< t_2$ for all $t_1\in I_1$, $t_2\in I_2$.
We then set $I_1\leq I_2$ whenever $I_1<I_2$ or $I_1=I_2$. The same notation will be used
for the corresponding characteristic functions $f=\chi_{I_{1}},\;g=\chi_{I_{2}}$, 
i.e. $f<g$ and $f\leq g$.
Note that, as in the monotone case, it is not possible to obtain a similar
formula for arbitrary functions, or even for characteristic functions with arbitrary supports.

\begin{Example}\label{p10}
{\rm
Let $f:=f_1 = f_2 = f_5 = f_6 = \chi_{[1,2)}$ and $g:=f_3 = f_4 = \chi_{[0,1)}$.
Besides, to simplify notation, we set $a^{\epsilon}(f_i) = a^{\epsilon}_i$ for $i=1,\ldots, 6$ and 
$\epsilon = 1,*$. Then
\begin{eqnarray}
\varphi(\omega(f_1) \omega(f_2) \ldots \omega(f_6) )
&=& 
\varphi( a^*_1 a_2 a^*_3 a_4 a^*_5 a_6)
+
\varphi(a^*_1 a^*_2 a_3 a_4 a^*_5 a_6)
\nonumber \\
&& 
+ 
\varphi(a^*_1 a_2 a^*_3 a^*_4 a_5 a_6)
+ 
\varphi(a^*_1 a^*_2 a_3 a^*_4 a_5 a_6)
\ \ \ \ \ \label{r36} \\ 
&& 
+
\varphi(a^*_1 a^*_2 a^*_3 a_4 a_5 a_6)
\nonumber \,,
\end{eqnarray}
since the other mixed moments certainly vanish. 
However, the second, third and fourth summands also give zero contribution
to (\ref{r36}) since in each of them a creation operator associated with $g$
is paired with an annihilation operator associated with $f$ or vice versa and these functions have
disjoint supports.
Therefore, it is enough to compute the contribution from the first and last
summands:
\begin{eqnarray*}
\varphi( a^*_1 a^*_2 a^*_3 a_4 a_5 a_6)
&=& 
\left\langle  a^*_1 a^*_2 a^*_3 
(g\oasterisk f\oasterisk f), \Omega \right\rangle =\ p\big\langle a^*_1 a^*_2 (f\oasterisk f), \Omega \big\rangle \\
&=& 
p\int_1^2 \int_1^2 w(t_1,t_2)dt_1dt_2
=\  \frac{pq + p}{2}\,.\\
\varphi(a^*_1 a_2 a^*_3 a_4 a^*_5 a_6) &=& 1.
\end{eqnarray*}
Thus
$$
\varphi(\omega(f_1) \omega(f_2) \ldots \omega(f_6))= \frac{pq + p + 2}{2}\,.
$$
}
\end{Example}

Before we find a formula for all moments,
let us introduce some additional notations (most of them are taken from [16]).

\begin{Definition}\label{d21}
{\rm
Let $f_1,f_2, \ldots, f_n \in\Theta$ have pairwise identical or disjoint supports.
We will say that $P=(P_1,\ldots,P_m)\in \mathcal{OP}_{n}$ is {\it adapted} to  
$(f_1,f_2,\ldots,f_n)$, which we denote $P\sim (f_1,f_2, \ldots,f_n)$, 
if and only if it satisfies two conditions:

1. $i,j\in P_k \ \Longrightarrow \ f_i=f_j\,,$

2. $i\in P_k$, $j\in P_l$ and $k<l \ \Longrightarrow \ f_i \leqslant f_j\,.$ 

\noindent
In turn, if $\pi=\{\pi_1,\ldots,\pi_m\}\in {\mathcal P}_{n}$, then we will 
say that $\pi$ is {\it adapted} to $(f_1,f_2,\ldots,f_n)$ if and only if its blocks
satisfy only the first condition, which we denote $\pi \sim (f_1,\ldots,f_n)$. }
\end{Definition}
\begin{Definition}
{\rm If $\pi\sim (f_1,f_2, \ldots , f_n)$, then the {\it support} of block $\pi_i$ 
is $\text{supp}\,\pi_i = \text{supp}\,f_j$, for any $j\in \pi_i$.
In turn, the {\it support} of the partition $\pi\in {\mathcal P}_{n}$ is
$$
\text{supp}\,\pi = \{(t_1,\ldots,t_m); t_k\in\text{supp}\,\pi_k,\, k=1,\ldots,m\}\,. %\text{supp}\,f_j\,,
$$
}
\end{Definition}

Finally, for $\pi \in \mathcal{NC}_{2k}^2$ and any $f_1,f_2, \ldots, f_{2k} \in\Theta$
we will also use a simplified notation
\begin{equation}\label{r41}
a_{\pi}(f_1,f_2, \ldots, f_{2k}) = a^{\epsilon_1}(f_1) a^{\epsilon_2}(f_2) \ldots a^{\epsilon_{2k}}(f_{2k})\,,
\end{equation}
where we understand that $\epsilon_i = *$  and $\epsilon_j = 1$ whenever $\{i,j\}$ is a block
of $\pi$ and $i<j$.

\begin{Lemma}\label{l13}
If $f_1,f_2, \ldots, f_{2n} \in\Theta$ have pairwise identical or disjoint 
supports and $\pi \in \mathcal{NC}_{2n}^2$ is not adapted to 
$(f_1,f_2, \ldots, f_{2n})$, then 
$\varphi\left( a_{\pi}(f_1,f_2, \ldots, f_{2n})\right)= 0$.
\end{Lemma}

\noindent
{\bf Proof.} 
Since $\pi$ is not adapted to $(f_1,\ldots,f_{2n})$, there exists a block
$\{i,j\}\in\pi$, such that $f_i,f_j$ have disjoint supports.
Assuming that $i<j$, we obtain
\begin{eqnarray*}
&&
\varphi\left( a_{\pi}(f_1,f_2, \ldots, f_{2n}) \right)\\
&=&
\left< a^{\epsilon_1}(f_1) \ldots a^{*}(f_i)  \ldots  a(f_j) \ldots a^{\epsilon_{2n}}(f_{2n}) \Omega, \Omega \right>\\[3pt]
&=& 
\left< a^{\epsilon_1}(f_1) \ldots a^{*}(f_i) \ldots a(f_j)  (g_1\oasterisk \ldots \oasterisk g_k), \Omega \right> \\[3pt]
&=& 
\left< a^{\epsilon_1}(f_1) \ldots a^{*}(f_i) \ldots a^{\epsilon_{j-1}}(f_{j-1}) 
(f_j\oasterisk g_1 \oasterisk \ldots \oasterisk g_k) , \Omega \right>\\[3pt]
&=& 
c \left< a^{\epsilon_1}(f_1) \ldots a^{*}(f_i) 
(f_j\oasterisk g_1 \oasterisk \ldots \oasterisk g_k),
\Omega \right> \\
&=& 
c \big< a^{\epsilon_1}(f_1) \ldots a^{\epsilon_{i-1}}(f_{i-1})
M(f_j,f_i) \,g_1 \oasterisk \ldots \oasterisk g_k, \Omega \big> \\
&=& 0\,,
\end{eqnarray*}
where we used the fact that $M(f_j,f_i)=0$ since $f_j$ and $f_i$ have
disjoint supports. $\hfill{\blacksquare}$

\vspace{10pt}
Let $\pi = \{ \pi_1,\ldots,\pi_k \}\sim (f_1, \ldots, f_{n})$, where
$f_{1},f_2, \ldots , f_n \in\Theta$ have pairwise identical or disjoint supports and
let $\sigma\in S_k$. Then the pair $(\pi,\sigma)$ can be identified
with an ordered partition, for which we can define
the set $\Delta_{(\pi,\sigma)}= \text{supp}\,\pi \cap \Delta_{\sigma}$, i.e.
\begin{equation}
\Delta_{(\pi,\sigma)} = \left\{ (t_1,\ldots,t_n)\in\mathbb{R}^n; t_{\sigma(1)} \leqslant \ldots \leqslant t_{\sigma(n)}\,, \ t_i\in \text{supp}\,\pi_i\,,\ i=1,\ldots,n \right\}.
\end{equation}
Then the following Proposition holds. 

\begin{Proposition}\label{l14}
Suppose $f_1,\ldots,f_{2n}\in\Theta$ have pairwise identical or disjoint supports
and let $\pi = \{ \pi_1,\ldots,\pi_n \}$ be a pair partition which is 
adapted to $(f_1,\ldots,f_{2n})$. Then, for any $\sigma \in S_n$ such that
$(\pi,\sigma)\nsim(f_1,\ldots,f_{2n})$ it holds that  $\Delta_{(\pi,\sigma)}=\emptyset$.
\end{Proposition}

\noindent
{\bf Proof.} We know that $\pi \sim (f_1,\ldots,f_{2n})$ and 
$(\pi,\sigma) \nsim (f_1,\ldots,f_{2n})$. Therefore, $(\pi,\sigma)$ deos not satisfy 
condition 2 of Definition \ref{d21}, i.e.
$$
\exists_{k,l\in\{1,\ldots,n\}} \ \forall_{i\in\pi_{\sigma(k)},j\in\pi_{\sigma(l)}} \ k<l \ \,\text{i}\, \ f_i>f_j.
$$ 
Let us suppose that $(t_1,\ldots,t_n) \in \Delta_{(\pi,\sigma)}$. Then it must hold that
\begin{equation}\label{r38}
t_{\sigma(k)} \leqslant t_{\sigma(l)}\,,
\end{equation}
since $k<l$. On the other hand, $f_i>f_j$ and $t_{\sigma(k)}\in \text{supp}\,\pi_{\sigma(k)} = \text{supp}\,f_i$ as well as $t_{\sigma(l)}\in \text{supp}\,\pi_{\sigma(l)} = \text{supp}\,f_j$, thus $t_{\sigma(k)} > t_{\sigma(l)}$, which contradicts (\ref{r38}). \hfill $\blacksquare$\\
\indent{\par}
In the sequel we will have to collect functions with the same supports in a suitable way. 
Suppose that $f_1,\ldots,f_{2n}\in\Theta$ have pairwise identical or disjoint supports and let
$g_1,\ldots,g_r\in\Theta$ be such that
\begin{equation}\label{r43}
\{f_1,\ldots,f_{2n}\} = \{g_1,\ldots,g_r\}\ \  \text{and}\ \ g_1 < \ldots < g_r\,
\end{equation}
and the corresponding (ordered) intervals by $I_1<I_2<\ldots <I_r$.
Then we can introduce numbers
\begin{equation}\label{r42}
b_i = \frac{1}{2}\big| \{ j\,; \text{supp\,}f_j = \text{supp\,}g^{(i)} \} \big|\ ,\ \ \ \ \ \ \ i=1,\ldots,r\,.
\end{equation}
Of course, $b_1 + b_2 + \ldots + b_r = n$. If there exists a pair partition which is
adapted to $(f_1,\ldots,f_{2n})$, then numbers $b_i$ are integers and the following
easy proposition holds.

\begin{Proposition}\label{l15}
If $f_1,\ldots,f_{2n}\in\Theta$ have pairwise identical or disjoint supports and
$(\pi,\sigma) = (\pi_{\sigma(1)},\pi_{\sigma(2)},\ldots,\pi_{\sigma(n)})\in \mathcal{ONC}_{2n}^2$ is adapted to $(f_1,\ldots,f_{2n})$, then
$$
\lambda(\Delta_{(\pi,\sigma)}) = 
\prod_{i=1}^{r}\frac{(\lambda(I_{i}))^{b_{i}}}{b_i!}
$$
where $\lambda$ denotes the (1-dimensional as well as n-dimensional) Lebesgue measure.
\end{Proposition}

\noindent
{\bf Proof.}
Since the volume of $\Delta_{(\pi,\sigma)}$ does not depend on $\sigma$, we have 
$$
\lambda(\Delta_{(\pi,\sigma)}) = \int_{\Delta_{(\pi,{\rm id})}} dt_1\ldots dt_n 
=
\prod_{i=1}^{r}\frac{(\lambda(I_{i}))^{b_{i}}}{b_i!}
$$
which gives the assertion.\hfill $\blacksquare$

\begin{Theorem}\label{t19}
If $f_1,f_2, \ldots, f_n \in\Theta$ have pairwise identical or disjoint supports, then
\begin{equation}\label{r40}
\varphi(\omega(f_1) \omega(f_2) \ldots \omega(f_n))
= 
\prod_{i=1}^{r}\frac{(\lambda(I_{i}))^{b_{i}}}{b_i!}
\sum_{\stackrel{P\in\mathcal{ONC}_{n}^2}{P \sim (f_1,\ldots, f_n)}} w(P)\,.
\end{equation}
\end{Theorem}

\noindent
{\bf Proof.} 
Of course, if $n$ is odd, the above moments vanish, which 
gives the assertion since (\ref{r40}) is a sum over the empty set.
Assume therefore that $n=2k$. First, observe that Lemma  \ref{l13} gives
\begin{eqnarray*}
\varphi(\omega(f_1) \omega(f_2) \ldots \omega(f_{2k}))
&=& 
\sum_{\pi\in \mathcal{NC}_{2k}^2} 
\varphi\left( a_{\pi}(f_1,f_2, \ldots, f_{2k})\right)\\
&=& \sum_{\stackrel{\pi\in\mathcal{NC}_{2k}^2}{\pi \sim (f_1,\ldots, f_{2k})}}
\varphi\left(a_{\pi}(f_1,f_2, \ldots, f_{2k}) \right)\,.
\end{eqnarray*}
If $\pi = \{ \pi_1,\ldots,\pi_k \} \in \mathcal{NC}_{2k}^2$, then 
to each number $i\in\{ 1,\ldots,k \}$ we can assign the number 
$l_i\in\{0,1, \ldots,k \}$, in such a way that block $\pi_{l_i}$ is a neighboring outer
block of $\pi_i$ whenever $\pi_i$ has any outer blocks, 
and otherwise we set $l_i = 0$ and $t_{0}=0$. Notice that
$$
\varphi\left(a_{\pi}(f_1,f_2, \ldots, f_{2k})\right)
= 
\int_{\text{supp\,} \pi} w(t_1,t_{l_1}) \ldots w(t_k,t_{l_k}) dt_1 \ldots dt_k\,.
$$
Therefore, we have
\begin{eqnarray*}
\varphi\left(
\omega(f_1) \omega(f_2) \ldots \omega(f_{2k})
\right)
&=& 
\sum_{\stackrel{\pi\in\mathcal{NC}_{2k}^2}{\pi \sim (f_1,\ldots, f_{2k})}} \int_{\text{supp\,} \pi} w(t_1,t_{l_1}) \ldots w(t_k,t_{l_k}) dt_1 \ldots dt_k \\
&=& 
\sum_{\stackrel{(\pi,\rho)\in\mathcal{ONC}_{2k}^2}{(\pi,\rho) \sim (f_1,\ldots, f_{2k})}} \int_{\Delta_{(\pi,\rho)}} w(t_1,t_{l_1}) \ldots w(t_k,t_{l_k}) dt_1 \ldots dt_k\\
&=& 
\prod_{i=1}^{r}\frac{(\lambda(I_{i}))^{b_{i}}}{b_i!}
\sum_{\stackrel{P\in\mathcal{ONC}_{2k}^2}{P \sim (f_1,\ldots, f_{2k})}} w(P)\,,
\end{eqnarray*}
using Propositions \ref{l14} and \ref{l15}.$\hfill{\blacksquare}$

\begin{Remark}
{\rm Note that there are two reasons why the mixed moments (\ref{r40})
depend on the supports of the $f_i$'s. The first one is the presence of the lenghts
of intervals $I_1, \ldots, I_r$ (in fact, in order that the RHS be not zero, 
each $f_{i}$ must appear an even number of times and therefore each lenght
can be replaced by an inner product).
The second one is the fact that the summation runs only over 
these $P$ which are adapted to $(f_1,f_2, \ldots, f_n)$.
For instance, if $f_1<f_2<\ldots <f_m$ and $g_1>g_2>\ldots >g_m$, we obtain
factorizations 
\begin{eqnarray*}
\varphi(\omega(f_1)\ldots \omega(f_m)\omega(f_m)\ldots \omega(f_1))&=&
q^{n-1}\varphi(\omega^{2}(f_1))\ldots \varphi (\omega^2(f_n))\\
\varphi(\omega(g_1)\ldots \omega(g_m)\omega(g_m)\ldots \omega(g_1))&=&
p^{n-1}\varphi(\omega^2(g_1))\ldots \varphi (\omega^2(g_n)).
\end{eqnarray*}
In particular, if $q=1$, the first property reflects the so-called 
`pyramidal factorization' of the mixed moments [12]. 
In fact, it is not hard to 
see that for $q=1$ we get pyramidal factorization for all mixed moments
of type $\varphi(c_1\ldots c_nd_n\ldots d_1)$, where $c_i,d_i$ are arbitrary 
elements of the unital algebra ${\mathcal A}_{i}=\langle 1,a(f_i),a^{*}(f_i) \rangle$,
$1\leq i \leq n$.
If $q\neq 1$, the first formula can be generalized to a 
`non-unital $q$-pyramidal factorization', in which $c_i,d_i\in \langle a(f_i),a^{*}(f_i)\rangle$.
A similar extension of the second formula is also possible.
In fact, one can show that for $q=1$ our processes satisfy all three conditions of the so-called 
`generalized Brownian motion' given in [6] (pyramidal factorization, stationarity and gaussianity).}
\end{Remark}
\begin{Example}
{\rm
In this example we will evaluate the sixth moment of $\omega(f)$ for 
$f=\chi_{[0,1)}$ by computing the contributions associated with all partitions 
$P\in \mathcal{ONC}^2_6$ corresponding to products of creation and annihilation operators.
These contributions will be compared with their counterparts obtained for the moments
of $t$-deformed operators studied in  [7]. We have
\begin{eqnarray*}
\varphi(\omega^6(f)) &=& 
\varphi( a^* a a^* a a^* a )
+ 
\varphi( a^* a^* a a a^* a ) 
+
\varphi(a^* a a^* a^* a a)\\
&&
+ 
\varphi(a^* a^* a a^* a a )  
+ 
\varphi(a^* a^* a^* a a a).
\end{eqnarray*}
These mixed moments of creation and annihilation operators are given 
in Table 1.
Analogous computations can be done for $t$-deformed creation and annihilation operators.
In Table 1 we compare the values of the summands in the above equation
for the $(p,q)$-interpolation and for the $t$-deformation.
The main observation is that, in general, the mixed moments of $(p,q)$-free creation
and annihilation operators are different than the corresponding mixed moments  
of $t$-deformed operators. In particular, for $p=0$ (arcsine law), 
the $t$-deformed moments do not reproduce the moments in the monotone case. }
\end{Example}
\begin{table}
\setlength{\tabcolsep}{3.0mm}
\renewcommand{\arraystretch}{1.5}

\vspace{15pt}
\begin{center}
\begin{tabular}{|c|c|c|c|c|c|}                                                        \hline
                     &  $(p,q)$-interpolation  & $t$-deformation &  $t\rightarrow (p+q)/2$  \\ \hline
$a^* a a^* a a^* a$  &      $1$            &      $1$         &       $1$                  \\ \hline
$a^* a^* a a a^* a$  &  $(p+q)/2$         &     $t$        &     $(p+q)/2$            \\ \hline
$a^* a a^* a^* a a$  &  $(p+q)/2$         &     $t$        &     $(p+q)/2$            \\ \hline
$a^* a^* a a^* a a$  &  $(p^2+pq+q^2)/3$     &    $t^2$       &     $(p^2+2pq+q^2)/4$       \\ \hline
$a^* a^* a^* a a a$  &  $(p^2+4pq+q^2)/6$    &    $t^2$       &     $(p^2+2pq+q^2)/4$       \\ \hline
                                                                                         \hline
$\sum$               &  $((p+q)^2+2p+2q+2)/2$    & $2t^2+2t+1$    &     $((p+q)^2+2p+2q+2)/2$       \\ \hline
\end{tabular} 
\end{center}
\caption{Comparison of mixed moments.}
\end{table}
Moreover, we can compare two combinatorial formulas for the moments.
It follows from [7] that the even moments of Kesten laws (1.4) satisfy the equation 
\begin{equation}\label{r51}
\mu_{2n} 
= 
\sum_{\pi\in \mathcal{NC}^2_{2n}} t^{{\rm in}(\pi)} 
= 
\sum_{k=0}^{n-1} {\mathcal D}(n,k) t^k 
\end{equation}
for $n\geq 1$, where $t=(p+q)/2$ and ${\rm in}(\pi)$ is the number of inner blocks in $\pi$ and the
${\mathcal D}(n,k)$ are the so-called {\it Delaney's numbers}, which give the numbers
of pair partitions in $\mathcal{NC}^2_{2n}$ which have exactly $k$ inner blocks.
They are given by the explicit formula
$$
{\mathcal D}(n,k) = {n+k-1 \choose k}  - {n+k-1 \choose k-1},
$$
where  
$k\in\{0,1,\ldots,n-1\}$, $n\geqslant 1$,
and $\left( \genfrac{}{}{0pt}{}{n}{-1} \right)=0$. 

\begin{Definition}\label{d23}
{\rm For $n\geqslant1$, the numbers of 
partitions from the set $\mathcal{ONC}_{2n}^2$ which have exactly $k$ disorders 
and $j$ orders will be called {\it generalized Euler's numbers} and will be denoted ${\mathcal E}(n,k,j)$, i.e.
$$
{\mathcal E}(n,k,j) = \left| \{ P\in\mathcal{ONC}_{2n}^2; e(P) = k\;{\rm and}\; e'(P)=j\} \right|
$$
}
where $0 \leq j,k \leq n-1$.
\end{Definition}

Using [7] and the results of this Section, we can 
find a relation between Delaney's numbers and generalized Euler's numbers.

\begin{Proposition}
For $n\geqslant1$ and $k\in\{0,1,\ldots,n-1\}$, it holds that
$$
{\mathcal E}(n,k,j) = 
\frac{n!}{2^{k+j}} 
{k+j \choose k} {\mathcal D}(n,k+j)\,.
$$
\end{Proposition}

\noindent
{\bf Proof.} Substituting $t=(p+q)/2$ in (\ref{r51}) and performing 
elementary algebraic calculations, we obtain
$$
\mu_{2n} = 
\sum_{k=0}^{n-1} \sum_{l=k}^{n-1} 
\frac{{\mathcal D}(n,l)}{2^l} 
{l \choose k} (p+q)^k\,,
$$
for $n\geqslant 1$. On the other hand, we know that
$$
\mu_{2n} = \frac{1}{n!} \sum_{P\in\mathcal{ONC}_{2n}^2} p^{e(P)} = 
\sum_{k,j=0}^{n-1} \frac{{\mathcal E}(n,k,j)}{n!} (p+q)^{k}
$$
Comparing the coefficients of these two polynomials, we get the desired
relation.\hfill$\blacksquare$

\section{Central limit theorem on discrete free Fock space}

In this section we will define discrete free creation and annihilation
operators on the free Fock space ${\mathcal F}(\mathcal{H})$ over a
Hilbert space $\mathcal{H}$ with a fixed orthonormal basis $\{e_i\}_{i=1}^{\infty}$ (for a discussion
of the discrete free Fock space, see [23]). 
Using them, we will formulate an elementary version of the central limit theorem 
for `$(p,q)$-independent' random variables. An abstract treatment of the notion of independence
involved here goes beyond the scope of this article and will be given in a separate paper.

Using the infinite matrix
$$
w_{i,j}= \left\{
\begin{array}{ll}
p & {\rm if}\;\; i<j\\
q & {\rm if}\;\; i>j\\
1 & {\rm if}\;\; i=j
\end{array}
\right.,
$$
we define a rescaled orthogonal basis in ${\mathcal F}(\mathcal{H})$,
$$
e_{i_1}\oasterisk e_{i_2}\oasterisk \ldots \oasterisk e_{i_n}
=
\sqrt{w_{i_1,i_2}w_{i_2,i_3}\ldots w_{i_{n-1},i_n}}
e_{i_1}\otimes e_{i_2}\otimes \ldots \otimes e_{i_n}
$$
where $i_1,i_2, \ldots , i_n\in {\mathbb N}$.
It is easy to see that 
$$
\langle 
e_{i_1}\oasterisk e_{i_2}\oasterisk \ldots \oasterisk e_{i_n},
e_{j_1}\oasterisk e_{j_2}\oasterisk \ldots \oasterisk e_{j_n}
\rangle
=
\delta_{i_1,j_1}\delta_{i_2,j_2}\ldots \delta_{i_n,j_n}
w_{i_1,i_2}w_{i_2,i_3}\ldots w_{i_{n-1},i_n}
$$

Using this basis, we define the creation operators
$A_i, i\in\mathbb{N}$ by equations
\begin{eqnarray*}
A_i (\Omega) &=& e_i \\
A_i \left(e_{i_1} \oasterisk e_{i_{2}}\ldots \oasterisk e_{i_n} \right)
&=& 
e_{i}\oasterisk e_{i_1}\oasterisk \ldots \oasterisk e_{i_n}
\end{eqnarray*}
with their adjoints, annihilation operators, 
acting as follows:
\begin{eqnarray*}
A_i^* (\Omega) &=& 0 \\
A_i^* (e_j) &=& \delta_{i,j}\Omega \\
A_i^* \left(e_{i_1} \oasterisk e_{i_2} \oasterisk \ldots \oasterisk e_{i_n} \right)
&=& w_{i_1,i_2}
\delta_{i,i_1} \, \, 
(e_{i_{2}} \oasterisk e_{i_3}\oasterisk \ldots \oasterisk e_{i_n} )
\end{eqnarray*}

\vspace{12pt}
We will study the asymptotic behavior of normalized sums
\begin{equation}\label{r60}
S_N = \frac{1}{\sqrt{N}} \sum_{i=1}^N \omega_i\,,
\end{equation}
where $\omega_i = A_i + A_i^*$ for $i\in\mathbb{N}$. 
It is easy to see that the position operators $\omega_i$ have mean zero and variance one
with respect to the vacuum state $\varphi$ on ${\mathcal F}({\mathcal H})$.
They will play the role of `independent' random variables in the central limit theorem.
In the propositions given below we state their properties which are rather standard in the central
limit context.

\begin{Proposition}\label{l20}
If, among the indices $i_1,i_2, \ldots, i_n \in {\mathbb N}$, there exists 
$i_{j}$, $j\in [n]$, such that $i_j \neq i_k$ for and $k\neq j$, then
$\varphi(\omega_{i_1} \ldots \omega_{i_n})=0$.
If, in turn, $(i_1,i_2, \ldots, i_{2n})$ is a sequence of indices associated with 
a partition $P\in\mathcal{OP}_{2n}^2$, then
$$
\varphi(\omega_{i_1} \ldots \omega_{i_{2n}})=
\begin{cases}
w(P)      &  \ \ \text{\rm if}\ P\in \mathcal{ONC}_{2n}^2\\
0         &  \ \ \text{\rm if}\ P\notin \mathcal{ONC}_{2n}^2
\end{cases}\,.
$$
\end{Proposition}
\noindent
{\bf Proof.} 
The first assertion is the usual singleton condition, which clearly holds in our case.
In second assertion is obvious if $P\notin \mathcal{ONC}_{2n}^2$ (it easily follows from
the definition of the $A_i$ and the $A^*_i$). 
Suppose that $P\in \mathcal{ONC}_{2n}^2$. 
Then there exists $r$ such that $i_r=i_{r+1}\neq i_{r+2}\neq \ldots \neq i_n $. Therefore
\begin{eqnarray*}
\varphi(\omega_{i_1} \ldots \omega_{i_n}) 
&=& 
\langle \omega_{i_1} \ldots \omega_{i_{r-1}} 
A^*_{i_r} A_{i_r} \omega_{i_{r+2}} \ldots \omega_{i_n} \Omega, \Omega \rangle \\
&=& 
w_{i_{r+1},i_{r+2}}
\langle \omega_{i_1} \ldots \omega_{i_{r-1}} \omega_{i_{r+2}} \ldots \omega_{i_n} \Omega , \Omega\rangle \\
&=& 
w_{i_{r+1},i_{r+2}}
\varphi( \omega_{i_1} \ldots \omega_{i_{r-1}} \omega_{i_{r+2}} \ldots \omega_{i_n} ),
\end{eqnarray*}
where $w_{i_{r+1},i_{r+2}}$ is equal to $p$ or $q$, depending on whether 
$i_{r+1} < i_{r+2}$ (then block $\{r,r+1\}$ forms a disorder with its neighboring
outer block) or $i_{r+1}>i_{r+2}$ (then block $\{r,r+1\}$ forms an order
with its neighboring outer block), respectively,
and otherwise $w_{i_{r+1},i_{r+2}}=1$ (then block $\{r,r+1\}$
does not have outer blocks).
The assertion follows from induction. \hfill$\blacksquare$

\begin{Theorem}\label{t14}
It holds that
$$
\lim_{N\rightarrow\infty} \varphi\left(S_N^{2n}\right) \ = 
\ \frac{1}{n!} \sum_{P\in \mathcal{ONC}_{2n}^2} w(P)
$$
and the odd moments vanish in the limit.
\end{Theorem}

\noindent
{\bf Proof.} 
The proof is standard since the mixed moments of the 
$\omega_{i}$ are invariant under order preserving injections, which gives
$$
\varphi\left(S_N^{2n}\right) = 
\frac{1}{N^n} \sum_{i_1,\ldots,i_{2n}}\varphi(\omega_{i_1}\omega_{i_2}\ldots \omega_{i_{2n}}) 
= 
\frac{1}{N^n} \sum_{r=1}^{2n} {N \choose r} 
\sum_{P\in \mathcal{OP}_{2n}(r)} \varphi(\omega_P)\,.
$$
where $\mathcal{OP}_{2n}(r)$ is the set of ordered partitions of the set $[n]$ which have 
$r$ blocks and $\varphi(\omega_{P})$ denotes 
$\varphi(\omega_{i_1}\omega_{i_2}\ldots \omega_{i_n})$ 
for any sequence $(i_1,i_2,\ldots,i_n)$ associated with $P$.
Using Proposition 4.1 and standard arguments, we obtain
the assertion.
\hfill$\blacksquare$\\

\section{Poisson processes}

In this Section we shall introduce processes of Poisson type, denoted $(\gamma_{t})_{t\geq 0}$ 
which correspond to the $(p,q)$-Brownian motion studied in the previous Section.
The moments of $\gamma_{t}$ in the vacuum state on 
${\mathcal F}({\mathbb R}_{+})$ 
are given by the combinatorial formula
$$
\varphi (\gamma_t^{n}) = 
\sum_{P\in\mathcal{ONC}_{n}} 
\frac{t^b(P)}{b(P)!} w(P)\,,
$$
where $t>0$ ($p,q$ are supressed in the notation). 
Again, as in the case of position processes, for $(p,q)=(1,1)$ 
we obtain the moments of the free Poisson process [22] given by 
$$
\varphi(\gamma_t^{n}) = \sum_{\pi\in\mathcal{NC}_{n}} t^{b(\pi)}
$$ 
where $b(\pi)$ denotes the number of blocks of $\pi$.
In turn, for $(p,q)=(0,1)$ we get the moments of the monotone 
Poisson process 
$$
\varphi(\gamma_t^{n}) 
= \sum_{P\in\mathcal{MON}_{n}} \frac{t^b(P)}{b(P)!},
$$ 
see [20]. Therefore, the process 
$(\gamma_{t})_{t\geq 0}$ plays the role
of a natural interpolation between these two processes.

To construct $\gamma_{t}$ we shall use the gauge operator (Definition 3.3)
associated with the characteristic function $\chi_{[0,t)}$, namely
$m_{t}=M(\chi_{[0,t)})$, which is given by the explicit formula
\begin{eqnarray*}
m_{t} \Omega  &=& 0\,,\\
m_{t}\, \left(f_1\oasterisk f_2\oasterisk \ldots \oasterisk f_n\right)
&=&  
(\chi_{[0,t)} f_1)\oasterisk f_2\oasterisk \ldots \oasterisk f_n\,.
\end{eqnarray*}
Such operators were used also in the free case in a similar context [22], where 
the Poisson process was defined as $l_{t} + l_{t}^* + l_{t}^{*}l_{t} +m_{t}$, 
where $l_{t}$ denotes the free creation operator associated with 
the function $\chi_{[0,t)}$.

An analogous form of the Poisson process will be adopted for the $(p,q)$-interpolation, 
namely
\begin{equation}
\gamma_{t} = a_{t} + a^*_{t} + a^*_{t}a_{t} + m_{t}\,
\end{equation}
and will be called the $(p,q)$-{\it Poisson process}.
From Proposition 3.3 it follows that 
$a^*_{t}a_{t}=n_{t}=M(\chi_{[0,t)}, \chi_{[0,t)})$ 
is another gauge operator given by the explicit formula
\begin{eqnarray*}
n_{t} \Omega  &=& t \Omega\,,\\
n_{t} \,\left( f_1\oasterisk f_2\oasterisk \ldots \oasterisk f_n\right)
&=&  (Wf_1)\oasterisk f_2\oasterisk \ldots \oasterisk f_n
\end{eqnarray*}
where $W(s) = \int_0^{t}w(u,s)du = \big((p-q)s + qt \big)$.

\begin{Example}\label{p21}
{\rm
Let us compute low order moments of the $(p,q)$-Poisson process.
\begin{eqnarray*}
\left< \gamma_t^1 \Omega, \Omega\right> &=& t\,, \\
\left< \gamma_t^2 \Omega, \Omega\right> &=& t + t^2 \,, \\
\left< \gamma_t^3 \Omega, \Omega\right> &=& t + \frac{p+q+4}{2}t^2 + t^3\,, \\
\left< \gamma_t^4 \Omega, \Omega\right> &=& t + \frac{3p+3q+6}{2}t^2 + \frac{p^2+pq+q^2+3p+3q}{3}t^3 +t^4\,, \\
\left< \gamma_t^5 \Omega, \Omega\right> &=& t + (3p+3q+4)t^2 + \frac{11p^2+11pq +11q^2 +24p +24q +36}{6}t^3 \\ 
& & \ \ +\ \frac{3p^3+3p^2q+3pq^2 +3q^2 +8p^2+8q^2 +18p+18q+48}{12} t^4 + t^5 
\end{eqnarray*}
For instance, we get 
\begin{eqnarray*}
\left< \gamma_t^4 \Omega, \Omega\right> &=& t\left< \gamma_t^3 \Omega, \Omega\right> + 
\left< \gamma_t^3 \big(\chi(t_1)\big), \Omega\right> \\
&=& t \left< \gamma_t^3 \Omega, \Omega\right> + \left< \gamma_t^2 \big(\chi(t_2)\chi(t_1) + (W(t_1)+1)\chi(t_1) +  t\Omega\big), \Omega\right> \\
&=& t \left< \gamma_t^3 \Omega, \Omega\right> + t \left< \gamma_t^2 \Omega, \Omega\right> \\
& & + \big< \gamma_t \big((W^2(t_1)+3W(t_1)+1)\chi(t_1) + (t + \frac{p+1}{2}t^2)\Omega \big), \Omega\big> \\
&=& t^2 + \frac{p+q+4}{2}t^3 + t^4 + t^2 + t^3 \\
& & + \int_0^{t}(W^2(t_1)+3W(t_1)+1)dt_1 + t^2 + \frac{p+q}{2}t^3 \\
&=& t + \frac{3p+3q+6}{2}t^2 + \frac{p^2+pq+q^2+3p+3q +9}{3}t^3 +t^4\,.
\end{eqnarray*}
where, for simplicity, we denote $\chi = \chi_{[0,t)}$.}
\end{Example}

Before we compute all moments of $\gamma_{t}$, 
we introduce some notations.
Let $\pi=\{\pi_1,\ldots,\pi_b\}\in\mathcal{NC}_n$ be a partition of $[n]$. 
Let us divide the set $[n]$ into the following disjoint subsets:
\begin{eqnarray*}
\mathfrak{a}_{\pi}   &=& \big\{ i\in [n]; \exists_{r\in\{1,\ldots,b\}} \ |\pi_r|>1,\ i=\max \pi_r \big\} \\
\mathfrak{a}^*_{\pi} &=& \big\{ i\in [n]; \exists_{r\in\{1,\ldots,b\}} \ |\pi_r|>1,\ i=\min \pi_r \big\} \\
\mathfrak{m}_{\pi}   &=& [n] \setminus (\mathfrak{n}_{\pi} \cup \mathfrak{a}_{\pi} \cup \mathfrak{a}^*_{\pi})\\
\mathfrak{n}_{\pi}   &=& \big\{ i\in [n]; \exists_{r\in\{1,\ldots,b\}} \ \pi_r = \{i\} \big\} \,.
\end{eqnarray*}
In other words, the sets  $\mathfrak{a}^*_{\pi}$ and $\mathfrak{a}_{\pi}$ consist of 
left and right legs of the blocks of $\pi$, respectively, 
$\mathfrak{m}_{\pi}$ corresponds to the `middle' legs of the blocks 
$\pi_1, \ldots, \pi_b$ and
the set $\mathfrak{n}_{\pi}$ corresponds to singletons in $\pi$.

Using these sets, we can assign to each 
$\pi=\{\pi_1,\ldots,\pi_b\}\in\mathcal{NC}_n$ an operator $c_{\pi}= c_1 \ldots c_n$, where
$$
c_i =
\begin{cases}
a_{t}   & \ \ {\rm if}\ \ i\in\mathfrak{a}_{\pi} \\
a^*_{t} & \ \ {\rm if}\ \ i\in\mathfrak{a}^*_{\pi} \\
m_{t}   & \ \ {\rm if}\ \ i\in\mathfrak{m}_{\pi} \\
n_{t}   & \ \ {\rm if}\ \ i\in\mathfrak{n}_{\pi} 
\end{cases}\,.
$$
Of course, the mapping $\pi \rightarrow c_{\pi}$ is one-to-one, but it is not onto.
However, the products of operators  
$a_{t}, a^*_{t}, n_{t}, m_{t}$, which do not correspond to 
any non-crossing partition $\pi$ will turn out irrelevant.
In Fig.4 we show a non-crossing partition $\pi \in \mathcal{ONC}_{10}(5)$ 
with the corresponding operator $c_{\pi}$. 
In turn, examples of products which do not correspond to any non-crossing partitions are
given by $m_{t}m_{t}m_{t}$, $a^*_{t}a_{t}a_{t}$, 
$a_{t}a^*_{t}a^*_{t}a_{t}$.

\begin{figure}
\unitlength=0.8mm
\special{em:linewidth 0.4pt}
\linethickness{0.4pt}
\begin{picture}(111.00,25.00)(100.00,35.00)

\put(165,45){$c_{\pi}=a^*_{t} a^*_{t} n_{t} m_{t} n_{t} a^*_{t} a_{t} a_{t} a_{t} n_{t}$}

\put(148,44){$\longrightarrow$}

\put(96.00,41.00){\line(0,1){4.00}}
\put(84.00,55.00){\line(1,0){48.00}}
\put(132.00,41.00){\line(0,1){14.00}}

\put(84.00,41.00){\line(0,1){14.00}}
\put(108.00,41.00){\line(0,1){8.00}}
\put(90.00,41.00){\line(0,1){8.00}}

\put(102.00,41.00){\line(0,1){8.00}}
\put(90.00,49.00){\line(1,0){36.00}}
\put(126.00,41.00){\line(0,1){8.00}}

\put(114.00,41.00){\line(0,1){4.00}}
\put(114.00,45.00){\line(1,0){6.00}}
\put(120.00,41.00){\line(0,1){4.00}}

\put(138.00,41.00){\line(0,1){4.00}}

\put(84.00,41.00){\circle*{1.00}}
\put(90.00,41.00){\circle*{1.00}}
\put(96.00,41.00){\circle*{1.00}}
\put(102.00,41.00){\circle*{1.00}}
\put(108.00,41.00){\circle*{1.00}}
\put(114.00,41.00){\circle*{1.00}}
\put(120.00,41.00){\circle*{1.00}}
\put(126.00,41.00){\circle*{1.00}}
\put(132.00,41.00){\circle*{1.00}}
\put(138.00,41.00){\circle*{1.00}}

\put(83.10,37.00) {\footnotesize 1}
\put(89.10,37.00) {\footnotesize 2}
\put(95.10,37.00) {\footnotesize 3}
\put(101.10,37.00) {\footnotesize 4}
\put(107.10,37.00) {\footnotesize 5}
\put(113.10,37.00) {\footnotesize 6}
\put(119.10,37.00) {\footnotesize 7}
\put(125.10,37.00) {\footnotesize 8}
\put(131.10,37.00) {\footnotesize 9}
\put(137.10,37.00) {\footnotesize 10}
%\put(100,28){$\pi$}

\end{picture}
\caption{Partition $\pi$ and the corresponding operator $c_{\pi}$.}
\end{figure}
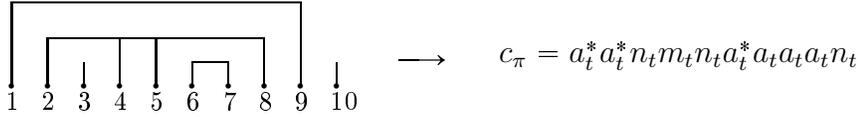

\begin{Lemma}\label{l21}
Let $c_1, \ldots, c_n \in \{ a_{t}, a^*_{t}, m_{t}, n_{t}\}$, 
be a sequence of operators, for which there exists
no partition $\pi\in\mathcal{NC}_n$ such that 
$c_{\pi}= c_1 \ldots c_n$. Then 
$\langle c_1 \ldots c_n \Omega, \Omega \rangle = 0$.
\end{Lemma}

\noindent
{\bf Proof.} 
Observe that if $c_{\pi} \neq c_1 \ldots c_n$ for all $\pi\in\mathcal{NC}_n$, 
then one of the following three cases must hold:

\noindent
1. The number of creation operators $a_{t}$ in the sequence 
$c_1, \ldots, c_n$ must be different from the number of annihilation operators
$a^*_{t}$. Then the assertion is certainly true.

\noindent
2. There exists $i\in [n]$ such that among $c_{i},c_{i+1}, \ldots, c_n$ there are 
more annihilation oparators $a^*_{t}$ than creation operators $a_{t}$. 
In that case, the assertion certainly holds.

\noindent
3. There exists $i\in [n]$ such that 
$c_{i}=m_{t}$ and in the sequence $c_{i+1}, \ldots, c_n$ there are as many annihilation 
operators $a^*_{t}$  as creation operators $a_{t}$. Then we have
\begin{eqnarray*}
\langle c_1 \ldots c_n \Omega, \Omega \rangle 
= 
\langle c_1 \ldots c_{i-1} m_{t} C \Omega, \Omega \rangle 
= 0\,,
\end{eqnarray*}
for a certain constant $C$. \hfill $\blacksquare$

\begin{Theorem}\label{t24}
The moments of the Poisson process $(\gamma_{t})_{t\geq 0}$ in the vacuum state are given by
$$
\varphi(\gamma_{t}^n)
= 
\sum_{P\in\mathcal{ONC}_{n}} \frac{t^b(P)}{b(P)!} w(P)
$$
for any $n\in\mathbb{N}$, and $\varphi(\gamma_{t}^{0})=1$.
\end{Theorem}

\noindent
{\bf Proof.} 
First, notice that from Lemma \ref{l21} it folows that 
$$
\varphi(\gamma_{t}^n)
= 
\sum_{\pi\in\mathcal{NC}_{n}} 
\left< c_{\pi}\Omega, \Omega\right>\,,
$$
and therefore it suffices to show that if $\pi\in\mathcal{NC}_{n}$ has 
$b(\pi)=b$ blocks, then it holds that
$$
\left< c_{\pi}\Omega, \Omega\right> 
= 
\frac{t^b}{b!} \sum_{\sigma\in S_b} w(\pi,\sigma)\,
$$
where $w(\pi, \sigma)=w(P)$  for $P=(\pi, \sigma)$.
For that purpose, let us construct a pair-partition $\pi''\in \mathcal{NC}^2_{2b}$ such that
$\left< c_{\pi}\Omega, \Omega\right> = \left< a_{\pi''}\Omega, \Omega\right>$ 
and  $e(\pi,\sigma)=e(\pi'',\sigma)$, $e'(\pi, \sigma)= e'(\pi'', \sigma)$
for any permutation $\sigma\in S_b$, where $a_{\pi''}$ is the
abbreviated notation for $a_{\pi''}(\chi_{[0,t)},\ldots, \chi_{[0,t)})$. 

Let us notice that in the product $c_{\pi}$ we can omit all occurences of 
$m_{t}$ since they correspond to the `middle' legs of blocks of $\pi$, and therefore
their action will not change the value of $\left< c_{\pi}\Omega, \Omega\right>$. 
In other words,
$$
\left< c_{\pi}\Omega, \Omega\right>
%= \langle \big( \prod^{\rightarrow}_{i\in\{1,\ldots,n\}} c_i \big) \Omega, \Omega \rangle_p 
= \langle \big( \prod^{\rightarrow}_{i\in\mathfrak{a}^*_{\pi} \cup \mathfrak{a}_{\pi} \cup \mathfrak{n}^*_{\pi}} c_i \big) \Omega, \Omega \rangle
= \left< c_{\pi'}\Omega, \Omega\right>\,,
$$
where $\pi'$ is a certain non-crossing partition with $b$ blocks, each 
consisting of one or two elements.
Substituting  the operator $a^*_{t} a_{t}$ for each $m_{t}$ in the product $c_{\pi'}$ corresponds 
to replacing singletons by two-element blocks $\{i,i+1\}$. 
Therefore, $c_{\pi'} = a_{\pi''}$ for a certain pair-partition $\pi''\in \mathcal{NC}^2_{2b}$.
Moreover, from the construction of $\pi''$ it follows that 
$e(\pi,\sigma)=e(\pi'',\sigma)$ and $e'(\pi, \sigma)=e'(\pi'', \sigma)$ 
for any permutation $\sigma\in S_b$. Now, from the proof of Theorem 3.1 we have
$$
\left< c_{\pi}\Omega, \Omega\right> 
= 
\left< a_{\pi''}\Omega, \Omega\right> 
= 
\frac{t^b}{b!} \sum_{\sigma\in S_b} w(\pi'',\sigma) 
= 
\frac{t^b}{b!} \sum_{\sigma\in S_b} w(\pi,\sigma)\,,
$$
which completes the proof. \hfill $\blacksquare$

\end{document}